\documentclass[11pt]{article}
\usepackage{amssymb,amsbsy,amsmath,amsfonts,amssymb,amscd}
\usepackage{latexsym}
\newtheorem{lemma}{Lemma}
\newtheorem{theorem}[lemma]{Theorem}
\newtheorem{corollary}[lemma]{Corollary}

\begin{document}

\def\C{{\mathbb C}}
\def\N{{\mathbb N}}
\def\Z{{\mathbb Z}}
\def\R{{\mathbb R}}
\def\PP{\cal P}
\def\p{\rho}
\def\phi{\varphi}
\def\ee{\epsilon}
\def\ll{\lambda}
\def\l{\lambda}
\def\a{\alpha}
\def\b{\beta}
\def\D{\Delta}
\def\g{\gamma}
\def\rk{\text{\rm rk}\,}
\def\dim{\text{\rm dim}\,}
\def\ker{\text{\rm ker}\,}
\def\square{\vrule height6pt width6pt depth 0pt}
\def\epsilon{\varepsilon}
\def\phi{\varphi}
\def\kappa{\varkappa}
\def\wz{\thinspace}
\def\proof{P\wz r\wz o\wz o\wz f.\hskip 6pt}
\def\quest#1{\hskip5pt {\scshape  Problem} {\rm #1}.\hskip 6pt}
\def\leq{\leqslant}
\def\geq{\geqslant}
\def\pd#1#2{\frac{\partial#1}{\partial#2}}
\def\limsup{\mathop{\overline{\hbox{\rm lim}\,}}}
\def\ug#1#2{\left\langle#1,#2\right\rangle}
\def\kk#1#2{{\k\langle#1,#2\rangle}}
\def\sv{\bf{ k} \langle X \rangle}
\def\k{k }
\def\lxr{\langle X \rangle}
\def\defin#1{\smallskip\noindent
{\scshape  Definition} {\bf #1}{\bf .}\hskip 8pt\sl}

\def\doubarr#1#2{\mathop{\hbox{$\vcenter{\offinterlineskip\halign
{\kern2pt\hfil##\hfil\kern2pt\cr \vrule height6pt width0pt
depth0pt\cr \smash{${\longleftarrow}\!\!{-}\!\!{-}$}\cr \vrule
height4pt width0pt depth0pt\cr
\smash{${\longleftarrow}\!\!{-}\!\!{-}$}\cr
}}$}}\limits_{#1}^{#2}}

\def\bull#1{\mathop{\bullet}\limits_{#1}}

\def\siglearr#1{\mathop{\hbox{${-}\!\!{-}\!\!{\longrightarrow}$}}\limits^{#1}}

\font\LINE=line10 scaled 1440 \font\Line=line10 \font\bIg=cmr10
scaled 1728 \font\biG=cmr10 scaled 1440 \font\LIne=line10
\def\lini{{\LINE\char"40}}
\def\li{{\Line\char"40}}
\def\lili{{\LIne\char"40}}
\def\pph{\vrule height4pt width0pt depth0pt}
\def\pha{{\phantom0}}
\def\sst{\scriptstyle}
\def\ddd{\displaystyle}
\def\>{\multispan}

\def\C{{\mathbb C}}
\def\N{{\mathbb N}}
\def\Z{{\mathbb Z}}
\def\R{{\mathbb R}}
\def\F{{\cal F}}
\def\U{{\cal U}}
\def\M{{\cal M}}
\def\Q{{\cal Q}}
\def\H{{\cal H}}
\def\E{{\cal E}}
\def\rr{{\cal R}}
\def\pp{{\cal P}}
\def\epsilon{\varepsilon}
\def\kappa{\varkappa}
\def\phi{\varphi}
\def\leq{\leqslant}
\def\geq{\geqslant}
\def\dim{\hbox{\rm dim}\,}
\def\ker{\hbox{\rm ker}\,}
\def\Cent{\hbox{\rm Cent}\,}
\def\ext{\hbox{\rm ext}\,}
\def\det{\hbox{\rm det}\,}
\def\deg{\hbox{\rm deg}\,}
\def\ssub#1#2{#1_{{}_{{\scriptstyle #2}}}}

\def\ramka#1#2#3{\hbox{$\vcenter{\offinterlineskip\halign
{\vrule\vrule\kern#2pt\hfil##\hfil\kern#2pt\vrule\vrule\cr
\noalign{\hrule}
\noalign{\hrule}
\vrule height #3pt depth0pt width0pt\cr
#1\cr
\vrule height #3pt depth0pt width0pt\cr
\noalign{\hrule}
\noalign{\hrule}}
}$}}

\def\jdots{\hbox{$\vcenter{\offinterlineskip\halign
{\hfil##\hfil\kern2pt&\hfil##\hfil\kern2pt&\hfil##\hfil\cr
&&.\cr
\noalign{\kern3pt}
&.&\cr
\noalign{\kern3pt}
.&&\cr
}}$}}

\vskip1cm

\title{Irreducible components of the Jordan varieties
}

\author{
 Natalia K. Iyudu}

\date{}

\maketitle

\small


\smallskip

\centerline{Department of Pure Mathematics, Queen's University
Belfast, Belfast BT7 1NN, U.K.}

\smallskip

\centerline{ {\bf e-mail:}  \,\, n.iyudu@qub.ac.uk}

\bigskip

\small

{\bf Abstract}

\bigskip

We announce here a number of results concerning representation theory of the algebra  $R=k\langle x,y\rangle/
(xy-yx-y^2)$, known as Jordan plane (or Jordan algebra).

We consider the question on 'classification' of finite-dimensional
modules over the Jordan algebra. Complete description of irreducible components of the
representation variety $mod (R,n)$, which we call a Jordan variety'  is given for any dimension $n$.
It is obtained on the basis of the stratification of this variety
related to the Jordan normal form of $Y$. Any irreducible component
of the representation variety contains only one stratum related to a
certain partition of $n$ and is the closure of this stratum. The
number of irreducible components therefore is equal to the number of
partitions of $n$.

As a preparation for the above result we describe the complete set of pairwise non-isomorphic irreducible
modules $S_{a}$ over the Jordan algebra, and the rule how they could be glued to
indecomposables. Namely, we show that  ${\rm
Ext}^1(S_{a},S_{b})=0$, if $a \neq b $.

We study then properties of the image algebras in the endomorphism
ring. Particularly, images of representations from the most important
stratum, corresponding to the full Jordan block $Y$. This stratum
turns out to be the only building block for the analogue of the
Krull-Remark-Schmidt decomposition theorem on the level of
irreducible components. Along this line we establish an analogue of
the Gerstenhaber--Taussky--Motzkin theorem on the dimension of
algebras generated by two commuting matrices. Another fact concerns
with the tame-wild question for those image algebras. We show that
all image algebras of $n$-dimensional representations are tame for $n \leq 4$
and wild for $n \geq 5$.


\normalsize

\vspace{15mm}

We consider here the quadratic algebra given by the presentation
$R=k\langle x,y\rangle/ (xy-yx-y^2)$. This algebra appears in
various different contexts in mathematics and physics. First of all,
it is a kind of quantum plane: one of the two Auslander regular
algebras of global dimension two, as it is mentioned in the Artin
and Shelter paper \cite{AS}. The other one is the usual quantum
plane $k\langle x,y\rangle/ (xy-qyx)$. It served as one of basic
examples for the foundation of noncommutative geometry (see
\cite{SM} and references therein). There were also studies of
deformations of $GL(2)$ analogous to $GL_q(2)$ with respect to $R$
in 80-90th  in Manin's 'Quantum group' \cite{Ma}, \cite{Ko},
 where this algebra appeared under the name
Jordan algebra.

This algebra is  also a simplest element in the class of  RIT
(relativistic internal time) algebras. The latter was defined and
investigated in \cite{a1}, \cite{ia}, \cite{wia},
\cite{iw}, \cite{cs}. The class of RIT algebras arises from a
modification of the Poincare algebra of the Lorenz group SO(3,1) by
means of introducing  the additional generator corresponding to
the relativistic internal time. The algebra $R$ above is a RIT
algebra of type (1,1). Our studies of this algebra are partially
reported in \cite{MP}.

 Let us mention
that  $R$ is a subalgebra of the first Weyl algebra $A_1$. The
latter has no finite dimensional representations, but $R$ turns out
to have quite a rich structure of them. Category of finite
dimensional modules over $R$ contains, for example, as a full
subcategory ${\rm mod} \, GP(n,2)$, where $GP(n,2)$ is the
Gelfand--Ponomarev algebra \cite{GP} with the nilpotency degrees of
variables $x$ and $y$,  $n$ and $2$ respectively. On the other hand
we show
that $R$ is residually finite
dimensional.

 We are interested here in
representations over an algebraically closed field $k$ of
characteristic $0$. Sometimes we just suppose $k= \Bbb C$, this will
be pointed out separately. We denote throughout the category of
all $R$-modules by $\rm{Mod }\,R$, the category of finite dimensional
$R$-modules by ${\rm mod}\, R$ and $\rho_n $ stands for an
$n$-dimensional representation of $R$.

The most important question on finite dimensional representations
one can ask, is to classify them. This question may be solved
directly by parametrization of isoclasses of indecomposable modules,
for example, for tame algebras. But when the algebra is wild the
problem turns to a description of orbits under $GL_n$ action by
simultaneous conjugation in the space $ mod(R,n)$ of
$n$-dimensional representations. One can think of the latter space
also as of  a variety    of tuples of $n\times n$ matrices
(corresponding to generators) satisfying the defining relations of
the algebra $R$.

It is commonly understood that the first step in the study of this
variety should be the description of its irreducible components.
This approach leads to such famous results of this kind as
Kashiwara-Saito \cite{K-S} description  of the irreducible
components of Lusztig's nilpotent variety via the crystal basis
\cite{Lus}.

We show that irreducible components of the representation space
$mod(R,n)$ of the Jordan algebra, which we call a Jordan space, could be completely described for
any dimension $n$.

The first key point for this description is the choice of a
stratification of the Jordan variety. We choose a stratum ${\cal
U}_{\cal P}$ corresponding to the partition ${\cal
P}=(n_1,\dots,n_k)$, $n=n_1+...+n_k, \,\, n_1 \ge...\ge n_k$,
consisting of the pairs of matrices $(X,Y)$, where $Y$ has the Jordan normal form defined by
${\cal P}$ and $X,Y$ satisfy the defining relation.

{\bf Theorem } {\it Any irreducible component $K_j$ of the
representation variety $mod(R,n)$ of the Jordan algebra contains only
one stratum $U_{\cal P}$ from the stratification related to the Jordan
normal form of $Y$, and is the closure of this stratum.

The number of irreducible components of the variety $mod(R,n)$ is
equal to the number of the partitions of $n$.

The variety $mod(R,n)$ is equidimensional, the dimension of components is $n^2$.
}

The importance of examples of algebras for which the irreducible
components of $ mod(R,n)$ could be described for each $n$ was
emphasized in \cite{CBSch} and it is mentioned there that known
cases are restricted to algebras of finite representation type
(i.e., there are only finitely many isomorphism classes of
indecomposable $R$ modules) and one example of infinite
representation type in \cite{GPS}.
There is some similarity between the algebra generated by the pair
of nilpotent matrices annihilating each other considered  in
\cite{GPS} and the Jordan algebra, but while in the case of
\cite{GPS} variables $x$ and $y$ act 'independently', there is much
more interaction in the case of Jordan algebra, which makes the
analysis in a sense more difficult.

We also can answer the question, in which irreducible components,
module in general position is indecomposable.

{\bf Corollary}  {\it Only the irreducible component
$K_{(n)}=\overline {{\cal U}_{(n)}}$ which is the closure of the
stratum corresponding to the trivial partition of $n$ (the full block
$Y$) contains an open dense subset consisting of indecomposable
modules}.

First, we derive some properties of
algebras which are images for representations of $R$ in the
endomorphism ring. We show that they are {\it basic} algebras, that
is their semisimple parts are direct sums of $r$ copies of the
field, where $r$ is the number of different eigenvalues of
$X=\rho_n(x)$. This allows to associate to any representation a {\it
quiver} of its image algebra, in a conventional way. This leads to a
rough classification of reps by these quivers. It turns out that
indecomposable modules have either a typical wild quiver with one
vertex and two loops or the quiver with one vertex and one loop. 

A simple, but important fact is that $Y=\rho_n(y)$ is nilpotent for
any $\rho_n \in mod \, R$. We actually prove slightly more general fact of linear algebra.
 
 \begin{lemma}
Let $X,Y$ be $n\times n$ matrices over an algebraically closed field
$k$ of characteristic zero. Assume that the commutator $Z=XY-YX$
commutes with $Y$. Then $Z$ is nilpotent.
\end{lemma}
 
 Note that this is not necessarily the
case when the characteristic of the basic field is not zero.

We describe the group of
automorphisms of the Jordan algebra. 

\begin{theorem}
 All automorphisms of $R=k\langle x,y|xy-yx=y^2\rangle$
are of the form $x\mapsto \alpha x+p(y)$, $y\mapsto \alpha y$, where
$\alpha\in k\setminus\{0\}$ and $p\in k[y]$ is a polynomial on $y$.
Hence the group of automorphisms isomorphic to a semidirect product
of an additive group of polynomials $k[y]$ and a multiplicative
group of the field  $k^*: \,\,{\rm Aut}(R) \simeq k[y]
\leftthreetimes k^*$.
\end{theorem}

On this bases we introduce later the
notion of auto-equivalence of modules: equivalence up to
automorphisms of the algebra $R$. 

We describe a
complete set of prime ideals of $R$ and point out which of them are
primitive.

\begin{theorem}
All prime ideals of $R=k\langle
x,y|xy-yx-y^2\rangle$ have a shape ${\rm id}(y)$ or ${\rm
id}(y,x-\alpha)$.
\end{theorem}

\begin{corollary}
\label{zmhsdfksj} 
The complete set of primitive
ideals in $R$ consists of the ideals ${\rm id}(y,x-\alpha)$,
$\alpha\in { k}$.
\end{corollary}

Then we develop main elements of representation theory for Jordan algebra.
We prove that all  {\it irreducible}
modules are one dimensional: $S_{\alpha}=(\alpha,0)$, and describe
all finite dimensional modules, subject to the Jordan normal form of
$Y$, in terms of $\cal B$-Toeplitz matrices. 

We study {\it
indecomposable} modules and derive the rule how one could glue
irreducibles together: ${\rm Ext}^1(S_{\a},S_{\b})=0$, if $\a \neq
\b $. These provides us with enough information to suggest a
stratification of the variety $mod(R,n)$ related to the Jordan
normal form of $Y$ and to prove results on irreducible components.

Results on the structure of representation variety for $R$ show an
exceptional role of the strata related to the full Jordan block
$Y$, since they turn out to be the only building blocks in the
analogue of the Krull-Remark-Schmidt decomposition theorem on the level of irreducible components.
Another
evidence of the special role of this stratum  is the following.
 We prove an analogue of the Gerstenhaber--Taussky--Motzkin
theorem \cite{Ger}, \cite{Gur} on the dimension of algebras
generated by two commuting matrices. We show that the dimension of
image algebras of representations of $R$ does not exceed $n(n+2)/4$
for even $n$ and $(n+1)^{2} /4$ for odd $n$ and this estimate is
attained in the stratum related the full Jordan block  $Y$.

Finally, we fulfill a
more detailed study of tame-wild questions for image algebras in the
strata related to the  full Jordan block  $Y$. From the above
results we see that any algebra which is an image of indecomposable
representation $\p_n$ is a {\it local complete} algebra. Hence after developing some technique,
which allows to write down defining relations in the image-algebra, we
apply Ringel's classification \cite{Ringloc} of local complete
algebras
and get that all image-algebras are {\it tame} for $n \leq 4$
and for $n \geq 5$ they are {\it wild}.

\vspace{15mm}

{\it Acknowledgments}

The work on this circle of questions was started
in 2003-2004
during my visit at the
Max-Planck-Intit\"ut f\"ur Mathematik in Bonn. Most of results apart from those on the irreducible components appeared in the MPI preprint \cite{MP}. I am thankful to this
institution for the support and hospitality and  to
many colleagues with whom I have been discussing ideas related to
this work.

\end{document}